# Ramanujan Primes and Bertrand's Postulate

## Jonathan Sondow

**1. INTRODUCTION: BERTRAND'S POSTULATE.** In a two-page paper [**8**; **9**, pp. 208-209] published one year before his death in 1920 at the age of 32, the Indian mathematical genius Srinivasa Ramanujan wrote:

> Landau in his *Handbuch* [**5**], pp. 89-92, gives a proof of a theorem the truth of which was conjectured by Bertrand: namely that there is at least one prime $p$ such that $x < p \leq 2x$, if $x \geq 1$. Landau's proof is substantially the same as that given by Tschebyschef. The following is a much simpler one.

Ramanujan then ``uses simple properties of the $\Gamma$-function'' (P. Ribenboim [**10**, p. 188]) to prove the theorem, which is known as *Bertrand's postulate* or *Tschebyschef's theorem*. (For an exposition of Ramanujan's proof, see Shapiro [**14**, Section 9.3B].)

**2. RAMANUJAN PRIMES.** In [**6**, p. 178] W. J. LeVeque explains that the theorem is called Bertrand's ``postulate''

> rather than ``conjecture'' because he took it as a working tool in his study of a problem in group theory. This must have seemed entirely safe, considering the actual density of primes in the tables. There is not merely one prime between 500,000 and 1,000,000, say, there are 36,960 of them!

This phenomenon is analyzed by Ramanujan at the end of his paper, where he proves the following extension of Bertrand's postulate. (I formulate it as a theorem and quote him.)

**Theorem 1 (Ramanujan).** ``*Let $\pi(x)$ denote the number of primes not exceeding $x$. Then ... $\pi(x) - \pi\left(\frac{1}{2}x\right) \geq 1, 2, 3, 4, 5, \ldots$, if $x \geq 2, 11, 17, 29, 41, \ldots$, respectively.*''

*Proof.* ``It follows from [inequalities established in the paper] that

$$\pi(x) - \pi\left(\tfrac{1}{2}x\right) > \frac{1}{\log x}\left(\tfrac{1}{6}x - 3\sqrt{x}\right), \text{ if } x > 300. \tag{1}$$

From this we easily deduce [the theorem].''  •

He means that from (1) the difference $\pi(x) - \pi\left(\frac{1}{2}x\right)$ tends to infinity with $x$. Notice that the case $\pi(x) - \pi\left(\frac{1}{2}x\right) \geq 1$, if $x \geq 2$, is Bertrand's postulate.

To describe the numbers $2, 11, 17, 29, 41, \ldots$, I coined the term ``Ramanujan primes'' [15, 16] in 2005. (Subsequently, I came across a completely different meaning of the term—see [7].)

**Definition 1.** For $n \geq 1$, the $n$th *Ramanujan prime* is the smallest positive integer $R_n$ with the property that if $x \geq R_n$, then $\pi(x) - \pi\left(\frac{1}{2}x\right) \geq n$.

Note that $R_n$ *is indeed a prime*, because by the minimality condition the functions $\pi(x) - \pi\left(\frac{1}{2}x\right)$ and, therefore, $\pi(x)$ must increase at $x = R_n$. Since they can increase by at most 1, *the equality* $\pi(R_n) - \pi\left(\frac{1}{2}R_n\right) = n$ *holds*.

**Example 1.** Bertrand's postulate is $R_1 = 2$. To compute $R_2$, set the quantity on the right side of inequality (1) equal to 1 and solve, obtaining $x = 392.39\ldots$. Since the quantity is an increasing function of $x$ in the range $x > 300$, and the left side of (1) is an integer and can change only at integers, it follows that $\pi(x) - \pi\left(\frac{1}{2}x\right) \geq 2$ if $x \geq 392$. This gives the bound $R_2 \leq 392$. Counting primes, we find that $\pi(n) - \pi\left(\frac{1}{2}n\right) \geq 2$ for $n = 391, 390, 389, \ldots, 11$, but $\pi(10) - \pi(5) = 1$. Thus $R_2 = 11$. In the same way, one gets $(R_3, R_4, R_5) = (17, 29, 41)$.

**3. BOUNDS FOR $R_n$.** The upper bound in the following result is much smaller than that derived from (1), and thus leads to a faster method of computing $R_n$.

**Theorem 2.** *The $n$th Ramanujan prime satisfies the inequalities*

$$2n \log 2n < R_n < 4n \log 4n \quad (n \geq 1).$$

*Furthermore, if $p_n$ denotes the $n$th prime, then $p_{2n} < R_n < p_{4n}$, for $n > 1$.*

*Proof.* We have $2\log 2 < R_1 = 2 < 4\log 4$. Now, by *Rosser's Theorem* [11], which asserts that $n \log n < p_n$, it suffices to show that if $n > 1$, then $p_{2n} < R_n < 4n \log 4n$.

To prove the first inequality, we first verify the cases $n = 2, 3, 4$. If $n \geq 5$, then in J. B. Rosser and L. Schoenfeld's [13] inequality $\pi(2x) < 2\pi(x)$, which holds for $x \geq 11$, we may take $x = \frac{1}{2}p_{2n}$, because $p_{10} = 29 > 22$. As $\pi(p_{2n}) = 2n$, we see that $\pi(p_{2n}) - \pi\left(\frac{1}{2}p_{2n}\right) < n$, and hence $p_{2n} < R_n$.

To prove the inequality $R_n < 4n \log 4n$, we first check it when $n < 4$. For $n \geq 4$, we use the lower bound [12]

$$\pi(2x) - \pi(x) > \frac{3}{5} \frac{x}{\log x} \quad (x \geq 20.5).$$

Notice that $x/\log x$ is increasing for $x > e$. Now take $x \geq 2n\log 4n$, and note that $x > 20.5$ because $n \geq 4$. Since by calculus

$$\frac{2y\log 4y}{\log(2y\log 4y)} > \frac{5}{3}y \quad (y \geq 4),$$

we get $\pi(2x) - \pi(x) > n$. Replacing $x$ with $\frac{1}{2}x$, we obtain that if $x \geq 4n\log 4n$, then $\pi(x) - \pi\left(\frac{1}{2}x\right) > n$. Therefore, $R_n < 4n\log 4n$. This completes the proof of the theorem. •

**Example 2.** We have $R_2 < 8\log 8 = 16.63\ldots$ (improving the bound $R_2 \leq 392$ from (1)).

**4. THE $n$TH RAMANUJAN PRIME IS ASYMPTOTIC TO THE $2n$TH PRIME.**
Using the *Prime Number Theorem* (PNT) [**5**, pp. 43-55; **6**, pp. 4-6; **10**, pp. 161-170], we improve the upper bound $R_n < 4n\log 4n$ by roughly a factor of 2, for $n$ large. We also show that the lower bound $p_{2n}$ is the true order of $R_n$.

**Theorem 3.** *For every $\varepsilon > 0$, there exists $n_0 = n_0(\varepsilon)$ such that*

$$R_n < (2 + \varepsilon)n\log n \quad (n \geq n_0). \tag{2}$$

*Moreover, we have the asymptotic formula $R_n \sim p_{2n}$ as $n \to \infty$, that is,*

$$\lim_{n\to\infty} \frac{R_n}{p_{2n}} = 1.$$

*Proof.* The PNT states that $\pi(x) \sim x/\log x$ as $x \to \infty$. It follows, using $\log(x/2) \sim \log x$, that given $\varepsilon > 0$, there exists $x_0 = x_0(\varepsilon) > 1$ such that

$$\left(\frac{1}{2} - \varepsilon\right)\frac{x}{\log x} < \pi(x) - \pi\left(\frac{1}{2}x\right) < \left(\frac{1}{2} + \varepsilon\right)\frac{x}{\log x} \quad (x \geq x_0).$$

From this we deduce that

$$(2 - \varepsilon)n\log n < R_n < (2 + \varepsilon)n\log n \quad (n \geq n_0),$$

for some $n_0 = n_0(\varepsilon)$. In particular, (2) holds, and $R_n \sim 2n\log n \sim 2n\log 2n$ as $n \to \infty$. The PNT implies $p_n \sim n\log n$ (see [**5**, p. 214] for a simple proof), and hence $R_n \sim p_{2n}$. •

**Example 3.** We have $R_{500}/p_{1000} = 8831/7919 = 1.115\ldots$.

**Example 4.** The number in LeVeque's illustration is

$$\pi(1,000,000) - \pi(500,000) = 36,960 = \left(\frac{1}{2} + 0.0106\ldots\right)\frac{1,000,000}{\log 1,000,000}.$$

**5. A CONJECTURAL UPPER BOUND.** The asymptotic formula $R_n \sim p_{2n}$ as $n \to \infty$ implies $R_n < p_{3n}$, for $n$ large. It also holds for $n \leq 1000$.

**Conjecture 1.** *The bound $R_n < p_{3n}$ holds true for all $n \geq 1$.*

(*Added in proof.* Shanta Laishram [4] has proved Conjecture 1. His proof uses P. Dusart's inequalities for Tschebyschef's function $\theta(x) := \sum_{p \leq x} \log p \leq \pi(x) \log x$, where $p$ is prime.)

Here is an unconditional result in the direction of Conjecture 1.

**Theorem 4.** *If $n \geq 1$, then $\pi(p_{3n}) - \pi\left(\frac{1}{2} p_{3n}\right) > n$.*

*Proof.* As $\pi(p_{3n}) = 3n$, we need to prove that $\pi\left(\frac{1}{2} p_{3n}\right) < 2n$. Invoking the inequality [12]

$$\pi(x) < \frac{5}{4} \frac{x}{\log x} \quad (x \geq 113.6),$$

we substitute $\frac{1}{2}x$ for $x$. Noting that $\log\left(\frac{1}{2}x\right) > \frac{15}{16} \log x$ when $x > 2^{16}$, we see that

$$\pi\left(\frac{1}{2}x\right) < \frac{2}{3} \frac{x}{\log x} \quad (x > 2^{2^4}).$$

Now take $x = p_{3n}$, with $n \geq 2181$ so that $x \geq p_{6543} = 2^{2^4} + 1 = F_4$ (the fourth and largest known Fermat prime [6, p. 3; 10, pp. 71-74]). From the inequality $p_k < k \log(k \log k)$ (valid [12] for $k \geq 6$) and Rosser's Theorem, we obtain $p_k < k \log p_k$, and hence

$$\pi\left(\frac{1}{2} p_{3n}\right) < \frac{2}{3} \frac{p_{3n}}{\log p_{3n}} < 2n \quad (n \geq 2181).$$

Checking that $\pi\left(\frac{1}{2} p_{3n}\right) < 2n$ also holds when $n < 2181$, the proof is complete. ●

**6. CONSECUTIVE RAMANUJAN PRIMES AMONG ALL PRIMES.** By Theorems 2 and 3, we have both $R_n > p_{2n}$ for $n > 1$, and $R_n \sim p_{2n}$ as $n \to \infty$. Thus, roughly speaking, *the probability of a randomly chosen prime being Ramanujan is slightly less than* $1/2$.

**Example 5.** We have $R_{500} = 8831 = p_{1100}$, so that $5/11$ of the first 1100 primes are Ramanujan primes.

Let us take a simplified model. According to S. Finch [**3**, p. 340], ``the expected length of the longest run of consecutive heads in a sequence of $n$ ideal coin tosses'' is approximately equal to

$$\frac{\log n + \gamma}{\log 2} - \frac{3}{2},$$

where $\gamma = 0.57\ldots$ is Euler's constant. For example, when $n = 1100$ the expected length is around 9.4.

By comparison, in a list of the first 1100 primes, the longest string of *consecutive Ramanujan primes* (respectively, consecutive non-Ramanujan primes) has length 13 (respectively, length 10). (It is a little surprising that the former is longer, because there are fewer Ramanujan primes than non-Ramanujan ones. On the other hand, while there is only one such string of length 13—namely, $(p_{384}, p_{385}, \ldots, p_{396}) = (R_{167}, R_{168}, \ldots, R_{179})$— the list contains three strings of non-Ramanujan primes of length 10.)

**Conjecture 2.** *In the sequence of prime numbers, there exist both arbitrarily long strings of consecutive Ramanujan primes, and arbitrarily long strings of consecutive non-Ramanujan primes.*

**7. TWIN RAMANUJAN PRIMES.** Recall that if $p$ and $p + 2$ are primes, they are called *twin primes* [**6**, p. 269; **10**, pp. 200-205]. We define *twin Ramanujan primes* simply as twin primes both of which are Ramanujan. Necessarily, they are of the form $R_n, R_{n+1}$, with $R_{n+1} = R_n + 2$. The smallest pair is $(R_{14}, R_{15}) = (149, 151)$.

If two primes are chosen at random, the probability that they are both Ramanujan is less than $1/2 \times 1/2 = 1/4$. But apparently the probability increases if they are twin primes.

**Example 6.** Among the first 1100 primes, there are 186 pairs of twin primes, and 70 pairs of twin Ramanujan primes. The ratio $70/186$ lies about midway between $1/4$ and $1/2$.

The following observation will help to explain why the ratio should be greater than may be expected a priori.

**Proposition 1.** *Let $p$ and $q = p + 2$ be twin primes greater than 5. Then*

$$\pi(q) - \pi\left(\tfrac{1}{2}q\right) = \pi(p) - \pi\left(\tfrac{1}{2}p\right) + 1. \tag{3}$$

*Proof.* It is easy to see that $(p, q) = (6k - 1, 6k + 1)$, for some integer $k > 1$. Since $3k$ is not prime, (3) follows. •

Now if *p* and *q* are *any* two primes with $p < q$, a necessary condition for them to be twin Ramanujan primes is evidently that equality (3) must hold. (The condition is not sufficient; in fact, neither of the twin primes 191 and 193 is Ramanujan.) Proposition 1 says that, in the case where *p* and *q* are twin primes greater than 5, the condition is automatically satisfied. Thus *two primes are more likely to be Ramanujan if they are twin primes*.

**Conjecture 3.** *For all* $x \geq 571$, *the inequality*

$$\frac{\#\{\text{pairs of twin Ramanujan primes} \leq x\}}{\#\{\text{pairs of twin primes} \leq x\}} > \frac{1}{4}$$

*holds. In particular, if there exist infinitely many twin primes, then there are also infinitely many twin Ramanujan primes.*

**8. ERDÖS'S PROOFS OF BERTRAND'S POSTULATE.** In his *Commentary on Ramanujan's Collected Works* [**9**, pp. 371-372], B. C. Berndt writes:

> P. Erdös, in his first published paper [**1**] in 1932, gave a proof of Bertrand's postulate ... which is quite similar to that of Ramanujan .... A short history of the origin of Erdös's proof is given in his paper *Ramanujan and I* [**2**, pp. 1-20].

In the latter, Erdös compares his proof to Ramanujan's: ``In fact the two proofs were very similar; my proof had perhaps the advantage of being more arithmetical." Ribenboim [**10**, p. 188] describes it as ``stressing divisibility properties of the middle binomial coefficient $\binom{2n}{n}$."

LeVeque [**6**, p. 178] says that a simpler proof than the one in [**1**] ``was found independently by Erdös and L. Kalmár in 1939, but was not published." (For the surprising reason, see [**2**].) An exposition of it is given in [**6**, Section 6.9].

At the end of [**1**], Erdös proves a version of Ramanujan's lower bound (1). Rewritten in the latter's notation, Erdös's inequality is

$$\pi(x) - \pi\left(\tfrac{1}{2}x\right) > \frac{\log 2}{60} \frac{x}{\log x}, \text{ if } x \geq 8000.$$

**ACKNOWLEDGEMENTS.** I thank both referees for several suggestions.

*209 West 97th Street, New York, NY 10025*
*jsondow@alumni.princeton.edu*